\documentclass[12pt]{amsart}

\usepackage{amsmath}
\usepackage{amssymb}
\usepackage{amscd}
\usepackage{indentfirst}
\usepackage[T1]{fontenc}
\usepackage[frenchb]{babel}

\theoremstyle{plain}
\newtheorem{thm}{Théorème}[section]
\newtheorem{lem}[thm]{Lemme}
\newtheorem{cor}[thm]{Corollaire}
\newtheorem{prop}[thm]{Proposition}
\theoremstyle{definition}
\newtheorem{dfn}[thm]{Définition}
\newtheorem{ntt}[thm]{Notations}
\theoremstyle{remark}
\newtheorem{rmq}[thm]{Remarque}

\DeclareMathOperator{\pic}{Pic}
\DeclareMathOperator{\ord}{ord}

\newcommand{\A}{\mathcal{A}}

\newcommand{\ez}{E^{\circ}}

\newcommand{\gm}{\mathbf{G}_{{\rm m}}}
\newcommand{\gmk}{\mathbf{G}_{{\rm m},K}}
\newcommand{\gmkk}{\mathbf{G}_{{\rm m},K'}}
\newcommand{\gmss}{\mathbf{G}_{{\rm m},S'}}
\newcommand{\gmt}{\mathbf{G}_{{\rm m},T}}

\newcommand{\Pc}{\mathfrak{P}}

\newcommand{\Homs}{{\rm Hom}}
\newcommand{\Hom}{\underline{\rm Hom}}
\newcommand{\Ext}{\underline{\rm Ext}^1}

\newcommand{\ext}{{\rm Ext}^1}

\DeclareMathOperator{\Biext}{BIEXT}
\DeclareMathOperator{\biext}{Biext^1}
\DeclareMathOperator{\Torsb}{TORSBIRIG}
\DeclareMathOperator{\Torsrig}{TORSRIG}

\begin{document}

\title{Invariants de classes : le cas semi-stable}

\author{Jean Gillibert}

\email{jean.gillibert@math.unicaen.fr}

\address{Université de Caen\\
Laboratoire de Mathématiques Nicolas Oresme (CNRS umr 6139)\\
BP 5186, 14032 CAEN cedex, FRANCE.}

\date{13 Janvier 2004}



\begin{abstract}
Nous construisons ici un analogue, pour un schéma en groupes semi-stable dont la fibre générique est une variété abélienne, du {\em class-invariant homomorphism} (introduit par M. J. Taylor dans le cadre des schémas abéliens), et nous en donnons une description géométrique. Puis nous généralisons un résultat de Taylor, Srivastav, Agboola et Pappas concernant le noyau de cet homomorphisme dans le cas d'une courbe elliptique.

\medskip

\noindent {\sc Abstract.} We define here an analogue, for a semi-stable group scheme whose generic fiber is an abelian variety, of M. J. Taylor's class-invariant homomorphism (defined for abelian schemes), and we give a geometric description of it. Then we extend a result of Taylor, Srivastav, Agboola and Pappas concerning the kernel of this homomorphism in the case of an elliptic curve.
\end{abstract}

\maketitle


\section{Introduction}

Soit $S$ un schéma, et soit $G\rightarrow S$ un $S$-schéma en groupes commutatif, fini et plat. Soit $G^D$ le dual de Cartier de $G$. Nous disposons d'un homomorphisme 
\[
\begin{CD}
\pi:H^1(S,G^D) @>>> \pic(G)\\
\end{CD}
\]
explicité en premier par Waterhouse (voir \cite{w}, Theorem $5$). Supposons que $S=Spec(R)$ soit affine. Alors on peut écrire $G=Spec(H)$, où $H$ est une $R$-algèbre de Hopf, et $G^D=Spec(H^*)$, où $H^*$ est l'algèbre duale de $H$. Si $X=Spec(C)$ est un $G^D$-torseur, alors $C$ est un $H^*$-comodule, donc un $H$-module. L'image de $X$ par $\pi$ est alors donnée par la classe de $C\otimes_H (H^*)^{-1}$ dans le groupe $\pic(H)=\pic(G)$.

Nous nous intéressons ici à un moyen de construction de $G^D$-torseurs dont l'image par $\pi$ est triviale, c'est-à-dire de torseurs dont la structure galoisienne est triviale.

Plus précisément, supposons que $S$ soit un schéma noethérien, excellent, intègre, normal, de point générique $\eta$. Soient $A$ et $A'$ deux $S$-schémas en groupes semi-stables, dont l'un des deux est à fibres connexes, et tels que $A_{\eta}$ et $A_{\eta}'$ soient deux variétés abéliennes duales l'une de l'autre. Supposons que $G$ soit un sous-groupe de $A$. On construit alors des $G^D$-torseurs grâce à la théorie de Kummer, puis on en déduit un homomorphisme
\[
\begin{CD}
\psi:A'(S) @>>> H^1(S,G^D) @>>> \pic(G)\,.\\
\end{CD}
\]
Le résultat principal de cet article est le suivant :

\begin{thm}
\label{jg}
Supposons que $A_{\eta}$ soit une courbe elliptique, et que l'ordre de $G$ soit premier à $6$. Alors $A'(S)_{\rm Tors}$ est contenu dans $\ker\psi$.
\end{thm}

Dans le cas particulier où $A$ est un $S$-schéma abélien de dimension 1 ({\em i.e.} une $S$-courbe elliptique), et où $G=A[m]$ (sous-groupe des points de $m$-torsion de $A$), ce résultat était déja connu. D'abord montré par Srivastav et Taylor dans \cite{st} pour une courbe elliptique à multiplication complexe sur un corps de nombres (en prenant $m$ égal à une puissance d'un nombre premier $\ell>3$), puis même sans l'hypothèse de multiplication complexe par Agboola dans \cite{a2}, il a été finalement prouvé par Pappas dans \cite{p1} (pour un schéma abélien de dimension 1 sur une base quelconque).

Dans le cas où $S$ est le spectre de l'anneau des entiers d'un corps de nombres $K$, le théorème \ref{jg} admet une interprétation arithmétique. Plus précisément, il implique que les anneaux d'entiers de certaines extensions de $K$, engendrées par des valeurs de fonctions elliptiques, sont libres en tant que modules sur l'algèbre de Hopf de $G$. Le lecteur trouvera dans \cite{cnt} de plus amples détails sur cette question, qui constitue la motivation initiale pour l'étude de ce problème.

Dans le cas où $A$ est une $S$-courbe elliptique et où l'ordre de $G$ est une puissance de $2$, Cassou-Noguès et Jehanne ont donné dans \cite{cnj} des exemples de non-annulation de $\psi$ sur les points de torsion. De plus, pour tout nombre premier $\ell$, Pappas a construit dans \cite{p1} une courbe affine lisse $S$ sur un corps fini et un $S$-schéma abélien $A$ de dimension $2$, tels que, si l'on pose $G=A[\ell]$, alors l'homomorphisme $\psi$ correspondant ne s'annule pas sur au moins un point de $\ell$-torsion.

Cependant la question reste ouverte pour un schéma abélien de dimension $\geq 2$ sur le spectre de l'anneau des entiers d'un corps de nombres. Nous nous contentons de signaler que $\psi$ s'annule toujours sur les points de torsion d'ordre premier à celui de $G$, quelle que soit la dimension de $A_{\eta}$ (voir la proposition \ref{oo}).

Le premier objectif de notre travail est la construction de $\psi$. Pour cela, nous avons dû adopter une approche différente de celle des auteurs précédents. Nous sommes amenés à utiliser (dans la section 2) une théorie de la dualité pour les schémas semi-stables qui s'énonce de la façon suivante : supposons à nouveau que $G\subseteq A$ soit un sous-groupe fini et plat de $A$ (pour des exemples, voir le paragraphe \ref{explicit}), et considérons la suite exacte
\[
\begin{CD}
0 @>>> G @>>> A @>\phi>> A/G @>>> 0\\
\end{CD}
\]
de faisceaux pour la topologie fppf sur $S$. En travaillant dans le petit site fppf des $S$-schémas plats (cf. paragraphe \ref{faisceaux}), nous vérifions l'exactitude de la suite
\[
\begin{CD}
0 @>>> G^D @>>> \Ext_S(A/G,\gm) @>\phi^*>> \Ext_S(A,\gm) @>>> 0\,.\\
\end{CD}
\]
Quand on restreint tous ces faisceaux à l'ouvert de bonne réduction de $A$, on retrouve un couple d'isogénies duales entre schémas abéliens.

Par application du foncteur des sections globales, on déduit de la suite précédente un morphisme cobord $\delta:\ext(A,\gm)\rightarrow H^1(S,G^D)$. On définit alors $\psi$ comme le composé des morphismes
\[
\begin{CD}
A'(S) @>\gamma>> \ext(A,\gm) @>\delta>> H^1(S,G^D) @>\pi>> \pic(G)\,,\\
\end{CD}
\]
la flèche $\gamma$ étant donnée par une biextension (voir le paragraphe \ref{bibi}). Si $A$ est un schéma abélien, et si $G=A[m]$, alors $A'$ est le schéma abélien dual de $A$, et $\psi$ est le {\em class-invariant homomorphism} de \cite{p1} (généralisant lui-même celui de Taylor \cite{t2}). Cette approche nous permet également, à travers quelques dévissages, de donner une preuve différente de la << description géométrique >> de $\psi$, généralisant celle de \cite{p1} (ce résultat apparaît en premier dans \cite{a1} dans le cas d'un schéma abélien sur le spectre de l'anneau des entiers d'un corps de nombres). Les précédents auteurs s'appuyaient sur des descriptions explicites de fibrés en droites sur les variétés abéliennes, tandis que notre construction permet de se ramener à l'étude de $\delta$. Ceci fait l'objet de la section 3.

Enfin, le but de la section 4 est de présenter une preuve du théorème \ref{jg}, que nous obtenons par des arguments analogues à ceux de Pappas \cite{p1}. Nous donnons ensuite deux applications de ce résultat.

\medskip

\noindent {\sc Remerciements}

Je tiens à remercier ici Laurent Moret-Bailly pour sa relecture du manuscrit et ses conseils. Je remercie également le Referee, dont les nombreuses suggestions m'ont permis d'améliorer et de clarifier l'exposition. Enfin je remercie vivement John Boxall pour son encadrement, ainsi que
Philippe Cassou-Noguès et Martin Taylor qui m'ont soutenu au cours de ce travail.


\section{Extension de l'isogénie duale}

Rappelons les notations qui seront en vigueur tout au long de cet article : $S$ est un schéma noethérien, excellent, intègre, normal, de point générique $\eta=Spec(K)$. Nous noterons $\gm$ le groupe multiplicatif sur $S$.

\begin{dfn}
\label{ssneron}
On dit qu'un $S$-schéma en groupes est semi-stable s'il est commutatif, lisse, séparé, et si les composantes neutres de ses fibres sont extensions de variétés abéliennes par des tores.

On dit qu'un $S$-schéma en groupes semi-stable $\A$ est néronien s'il vérifie la propriété universelle suivante : pour tout $S$-schéma lisse $Y$ et tout morphisme $\alpha_K:Y_{\eta}\rightarrow \A_{\eta}$, il existe un unique morphisme $\alpha:Y\rightarrow \A$ prolongeant $\alpha_K$.
\end{dfn}

On fixe une fois pour toutes un $S$-schéma en groupes semi-stable $A$ dont la fibre générique $A_{\eta}$ est une variété abélienne. On notera $U\subseteq S$ l'ouvert de bonne réduction de $A$, de sorte que $A_U$ est un $U$-schéma abélien.

\subsection{Faisceaux sur le << petit site fppf >>}
\label{faisceaux}
Soit $j:U\rightarrow S$ l'inclusion. On munit la catégorie des schémas plats sur $S$ (resp. sur $U$) d'une structure de site pour la topologie fppf (on obtient ce qu'on appelle un << petit site fppf >>, la catégorie sous-jacente étant la catégorie des $S$-schémas plats). Le choix de ce petit site nous permettra de démontrer le lemme \ref{null}.

Soit $j^*$ le foncteur << image inverse >>, qui à un faisceau sur $S$ associe sa restriction à l'ouvert $U$. Comme $j:U\rightarrow S$ est un objet du petit site fppf sur $S$, $j^*$ est un << foncteur de localisation >> (voir \cite{gro4}, exposé IV, paragraphes 5.1 à 5.4). En particulier, l'image par $j^*$ d'un faisceau représentable (disons par un $S$-schéma plat $X$) est représenté par le $U$-schéma $X_U:=X\times_S U$.

D'autre part, si $F_1$ et $F_2$ sont deux faisceaux abéliens sur $S$, la flèche canonique
$$j^*(\Hom_S(F_1,F_2))\rightarrow \Hom_U(j^*F_1,j^*F_2)$$
est un isomorphisme (voir \cite{gro4}, exposé IV, prop. 12.3, b), p. 502). De plus, $j^*$ est exact et admet un adjoint à gauche $j_{!}$ exact. Par suite, $j^*$ envoie les injectifs sur des injectifs (voir \cite{gro4}, exposé V, paragraphe 2.2). Nous dérivons alors des deux côtés de la flèche et obtenons un isomorphisme $j^*(\Ext_S(F_1,F_2))\simeq \Ext_U(j^*F_1,j^*F_2)$. En particulier, soit $X$ un $S$-schéma en groupes plat tel que $X_U$ soit un $U$-schéma abélien, alors nous obtenons un isomorphisme
\begin{equation}
\label{iminv}
j^*(\Ext_S(X,\gm))\simeq X_U^t
\end{equation}
où $X_U^t=\underline\pic_{X_U/U}^0$ est le schéma abélien dual de $X_U$ (voir \cite{fc}, chap. I, \S\/1). On se sert ici du fait que $\Ext_U(X_U,\gm)$ est isomorphe à $X_U^t$ (voir \cite{gro7}, exposé VII, 2.9.5 et 2.9.6).

\begin{rmq}
Soit $H$ un $S$-schéma en groupes commutatif affine, plat et localement de type fini (donc localement de présentation finie, $S$ étant noethérien). Alors les $H$-torseurs pour la topologie fpqc sont représentables. Un argument de descente fidèlement plate (voir \cite{ega4}, prop. 2.7.1) montre que les $H$-torseurs fpqc sont des torseurs fppf, d'où l'égalité $H^1_{\rm fpqc}(S,H)=H^1_{\rm fppf}(S,H)$. Un raisonnement analogue montre que, si $H'$ est un autre $S$-schéma en groupes commutatif, alors $\ext_{\rm fpqc}(H',H)=\ext_{\rm fppf}(H',H)$.
\end{rmq}

\subsection{Isogénies duales}
\label{isodu}
Soit $G_{\eta}$ un sous-groupe algébrique fini de $A_{\eta}$, et soit $B_{\eta}$ la variété abélienne quotient $A_{\eta}/G_{\eta}$, nous obtenons une isogénie $A_{\eta}\rightarrow B_{\eta}$ dont le noyau est égal à $G_{\eta}$. Alors le noyau de l'isogénie duale $B_{\eta}^t\rightarrow A_{\eta}^t$ s'identifie au dual de Cartier $G_{\eta}^D$ de $G_{\eta}$.

Ce résultat, classique sur un corps, s'étend au cas des schémas abéliens, pour lesquels on dispose d'une notion de schéma dual (voir \cite{fc}, chap. I, \S\/1).

On fixe à présent un sous-$S$-schéma en groupes fini et plat $G$ de $A$. Nous noterons $B$ le faisceau quotient $A/G$ pour la topologie fppf sur $S$.

\begin{rmq}
\label{dimsun}
Si $S$ est de dimension $\leq 1$, alors $B$ est représentable par un $S$-schéma en groupes, également semi-stable (voir \cite{anan}, chap. IV, théorème 4.C).
\end{rmq}

Nous avons (par définition) une suite exacte
\[
\begin{CD}
0 @>>> G @>>> A @>\phi>> B @>>> 0\\
\end{CD}
\]
de faisceaux abéliens (sur le petit site fppf de $S$), prolongeant la suite exacte
\[
\begin{CD}
0 @>>> G_U @>>> A_U @>\phi_U>> j^*B @>>> 0\\
\end{CD}
\]
dans laquelle $A_U$ est un $U$-schéma abélien. On en déduit que $j^*B$ est le quotient  $A_U/G_U$, donc est représentable par un $U$-schéma abélien, que nous noterons $B_U$.

Nous obtenons alors, en appliquant le foncteur $\Hom_S(-,\gm)$ à la première suite, une (longue) suite exacte de cohomologie
\begin{align*}
\Hom_S(A,\gm) \rightarrow \Hom_S(G,\gm) \rightarrow \Ext_S(B,\gm) \rightarrow \Ext_S(A,\gm) \rightarrow \Ext_S(G,\gm)
\end{align*}
dont les termes sont des faisceaux abéliens sur $S$.

D'autre part, $G$ étant fini et plat sur $S$, le faisceau $\Hom_S(G,\gm)$ est représentable par $G^D$ (le dual de Cartier de $G$). Pour les mêmes raisons, le faisceau $\Ext_S(G,\gm)$ est nul (voir \cite{gro7}, exposé VIII, 3.3.1). Enfin nous avons le lemme qui suit :

\begin{lem}
\label{null}
Soit $X$ un $S$-schéma en groupes plat, dont la fibre générique $X_K$ est une variété abélienne. Alors $\Hom_S(X,\gm)$ est nul.
\end{lem}

\begin{proof}
Nous devons montrer, pour tout schéma $S'\rightarrow S$ plat, la trivialité
du groupe $\Homs_{S'-gr}(X_{S'},\gmss)$. Pour
cela, il suffit de se limiter au cas où $S$ et $S$' sont affines. Soient donc
$S=Spec(R)$ et $S'=Spec(R')$ où $R'$ est une $R$-algèbre plate. Cette hypothèse de
platitude permet d'affirmer que
$\mathcal{O}_{X_{S'}}(X_{S'})= \mathcal{O}_X(X)\otimes_R R'$. On en déduit que
$\mathcal{O}_{X_{S'}}(X_{S'})$ est $R'$-plat, sachant que
$\mathcal{O}_X(X)$ est $R$-plat. D'autre part, soit
$K'=R'\otimes_R K$, alors $K'$ est une $K$-algèbre, et la flèche $R'\rightarrow K'$ est injective
par $R$-platitude de $R'$. Finalement, on trouve que la flèche
$$\mathcal{O}_{X_{S'}}(X_{S'})\rightarrow
\mathcal{O}_{X_{S'}}(X_{S'})\otimes_{R'} K'=
\mathcal{O}_{X_{K'}}(X_{K'})$$
est un morphisme injectif d'anneaux. Sa restriction $\mathcal{O}_{X_{S'}}(X_{S'})^{\times}\rightarrow
\mathcal{O}_{X_{K'}}(X_{K'})^{\times}$ est donc injective. En d'autres termes, l'application
$$\Homs_{S'}(X_{S'},\gmss)\rightarrow \Homs_{K'}(X_{K'},\gmkk)$$
obtenue par changement de base $Spec(K')\rightarrow S'$, est
injective. On en déduit que l'application
$\Homs_{S'-gr}(X_{S'},\gmss)\rightarrow \Homs_{K'-gr}(X_{K'},\gmkk)$
est également injective, en ne considérant que les morphismes de groupes.

Montrons à présent la trivialité de $\Homs_{K'-gr}(X_{K'},\gmkk)$. Soit $f:X_{K'}\rightarrow \gmkk$ un morphisme de $K'$-schémas en
groupes, $f$ est entièrement déterminé par la donnée du morphisme de
$K'$-algèbres de Hopf
$f^{\sharp}:K'[x,x^{-1}]\rightarrow
\mathcal{O}_{X_{K'}}(X_{K'})$. D'autre part,
$\mathcal{O}_{X_{K'}}(X_{K'})\simeq
\mathcal{O}_{X_{K}}(X_{K})\otimes_K K'$ et
$\mathcal{O}_{X_{K}}(X_{K})\simeq K$ (car $X_K$ est propre et
géométriquement intègre), donc
$\mathcal{O}_{X_{K'}}(X_{K'})\simeq K'$ en tant que $K'$-algèbres de
Hopf. La counité étant l'unique morphisme de $K'$-algèbres de Hopf
$K'[x,x^{-1}]\rightarrow K'$, ceci
prouve que $f^{\sharp}$ se factorise par la counité. Au final, $f$ est
le morphisme trivial. 
\end{proof}

En appliquant le lemme \ref{null} au schéma $A$, nous obtenons une suite exacte
\begin{equation}
\label{exsq}
\begin{CD}
0 @>>> G^D @>>> \Ext_S(B,\gm) @>\phi^*>> \Ext_S(A,\gm) @>>> 0\,.\\
\end{CD}
\end{equation}
D'après ce qui précède (voir (\ref{iminv})), son image par le foncteur $j^*$ est la suite
\[
\begin{CD}
0 @>>> G_U^D @>>> B_U^t @>\phi_U^t>> A_U^t @>>> 0\,.\\
\end{CD}
\]
Cette dernière admet donc un << prolongement >> sur $S$ en termes de faisceaux.

\subsection{Liens avec la théorie des biextensions}
\label{bibi}
Soit $A_{\eta}^t$ la variété abélienne duale de $A_{\eta}$, et soit $A'$ un $S$-schéma en groupes semi-stable prolongeant $A_{\eta}^t$. On sait qu'il existe un unique tel prolongement à fibres connexes (voir \cite{mb}, chap. IV, th. 7.1 (i)).

L'ouvert de bonne réduction de $A'$ est le même que celui de $A$, donc $(A')_U$ coïncide avec le schéma abélien dual $A_U^t$ de $A_U$.

Nous voudrions établir une dualité entre $A$ et $A'$ prolongeant celle qui existe déjà sur les fibres génériques. On sait que la dualité entre $A_{\eta}$ et $A_{\eta}^t$ découle de l'existence d'un fibré en droites $\mathcal{P}_{\eta}$ sur $A_{\eta}\times_K A_{\eta}^t$, que l'on appelle fibré de Poincaré. Cependant nous allons envisager ici la dualité dans un cadre plus général à l'aide de la notion de biextension, imaginée par Mumford dans \cite{mu}. Nous introduisons brièvement ici cet outil. Pour les détails nous renvoyons le lecteur à l'exposé de Grothendieck (\cite{gro7}, exposé VII). Mentionnons également Milne (\cite{mil}, Appendix C).

\begin{dfn}
Soient $P$ et $Q$ deux $S$-schémas en groupes commutatifs. On peut considérer sur $P\times_S Q$, en plus de sa structure naturelle de $S$-schéma en groupes, sa structure de $P$-schéma en groupes et celle de $Q$-schéma en groupes. Ces deux schémas en groupes seront notés respectivement $Q_P$ et $P_Q$.
Une {\em biextension} de $(P,Q)$ par $\gm$ est un $\gm$-torseur $Y$ sur $P\times_S Q$, muni de lois de composition partielles, qui fassent du $P$-schéma $Y_P$ une extension de $Q_P$ par $\mathbf{G}_{{\rm m},P}$, et du $Q$-schéma $Y_Q$ une extension de $P_Q$ par $\mathbf{G}_{{\rm m},Q}$. En outre, nous demandons que ces deux structures d'extensions soient compatibles en un certain sens (voir \cite{gro7}, exposé VII, \S\/2 ou \cite{mil}, Appendix C).
\end{dfn}

On note $\Biext(P,Q;\gm)$ la catégorie des biextensions de $(P,Q)$ par $\gm$.
Soit $Y$ une telle biextension, alors $Y$ définit un morphisme de groupes
$$Q(S)\longrightarrow \ext(P,\gm)$$
que nous explicitons : soit $f:S\rightarrow Q$ une section de $Q$, alors $({\rm id}_P \times f)^*(Y)$ est une biextension de $(P,e)$ par $\gm$ (où $e$ est le $S$-groupe trivial), c'est-à-dire une extension de $P$ par $\gm$.

D'autre part, on note $\biext(P,Q;\gm)$ le groupe constitué par l'ensemble des classes d'isomorphie de biextensions de $(P,Q)$ par $\gm$. On dispose d'un homomorphisme de groupes
$$t:\biext(P,Q;\gm)\rightarrow \pic(P\times_S Q)$$
qui à toute biextension associe son $\gm$-torseur sous-jacent. On constate que $t$ est une transformation naturelle entre (bi)foncteurs.

Retournons à notre question de dualité. Grâce au théorème du carré, on peut munir le fibré de Poincaré $\mathcal{P}_{\eta}$ d'une unique structure de biextension, que l'on appelle la biextension de Weil, et que l'on note $W_{\eta}$. Le problème se reformule alors de la façon suivante : peut-on prolonger la biextension $W_{\eta}\in \Biext(A_{\eta}, A_{\eta}^t; \gmk)$ en une biextension qui vive dans $\Biext(A,A'; \gm)$ ? Toujours d'après Grothendieck, il existe une obstruction à ce prolongement (voir \cite{gro7}, exposé VIII, théorème 7.1, b)), qui s'incarne sous la forme d'un accouplement (dit de monodromie) entre les groupes de composantes de $A$ et $A'$. Si $A$ ou $A'$ est à fibres connexes, cette obstruction disparaît :

\begin{prop}
\label{wo}
Supposons que $A$ ou $A'$ soit à fibres connexes. Alors il existe une unique biextension $W$ de $(A,A')$ par $\gm$ prolongeant la biextension de Weil $W_{\eta}$ sur $(A_{\eta},A_{\eta}^t)$.
\end{prop}

\begin{proof}
Considérons le foncteur de restriction à la fibre générique
$$\Biext(A,A'; \gm)\longrightarrow \Biext(A_{\eta}, A_{\eta}^t; \gmk)\,.$$
D'après (\cite{mb}, chap. II, th. 3.6), ce foncteur est une équivalence de catégories. On en déduit le résultat voulu.
\end{proof}

Nous supposons à partir d'ici, et dans toute la suite de ce travail, que $A$ ou $A'$ est à fibres connexes. Nous noterons $W$ la biextension dont l'existence est assurée par la proposition \ref{wo} ainsi que $\gamma:A'(S)\rightarrow \ext(A,\gm)$ le morphisme défini par $W$.

\begin{rmq}
\label{gamma}
Dans le cas où $A$ est un schéma abélien, alors $A'$ est le dual de $A$ et $\gamma$ est un isomorphisme. Dans le cas général, considérons le diagramme commutatif suivant :
\[
\begin{CD}
A'(S) @>\gamma>> \ext(A,\gm) \\
@VVV @VVV \\
A_{\eta}^t(K) @>>> \ext(A_{\eta},\gmk)\,. \\
\end{CD}
\]
D'après (\cite{gro7}, exposé VIII, théorème 7.1 et remarque 7.2), la flèche de droite est injective. Il est clair que la flèche de gauche est injective (car $A'$ est séparé sur $S$) et que celle du bas est bijective. Ainsi $\gamma$ est toujours injective. Si de plus $A$ est à fibres connexes et $A'$ est néronien (déf. \ref{ssneron}), alors les flèches verticales sont bijectives (celle de gauche d'après la propriété universelle de Néron, celle de droite d'après {\em loc. cit.}), donc $\gamma$ est un isomorphisme.
\end{rmq}


\section{L'homomorphisme $\psi$ et sa géométrie}

Reprenons les notations et hypothèses du début de la section 2. En appliquant à la suite exacte (\ref{exsq}) le foncteur des sections sur $S$, nous obtenons une suite exacte de cohomologie
\[
\begin{CD}
\cdots @>>> \Ext_S(A,\gm)(S) @>\delta>> H^1(S,G^D) @>>> \cdots\\
\end{CD}
\]
Ce cobord $\delta$ va jouer un rôle essentiel dans la définition de notre homomorphisme. Tout d'abord, il convient de mieux connaître les objets en jeu, ce qui va nous permettre de donner une autre description (dite << géométrique >>) de $\pi\circ\delta$.

\subsection{Description du cobord}
Commençons par énoncer un lemme de comparaison entre Ext locaux et globaux.

\begin{lem}
\label{ext}
Soit $X$ un $S$-schéma en groupes plat, dont la fibre générique $X_K$ est une variété abélienne. Alors le faisceau $\Ext_S(X,\gm)$ est canoniquement isomorphe au faisceau
$T\mapsto\ext(X_T,\gmt)$.
\end{lem}

\begin{proof}
Soient $F$ et $F'$ deux faisceaux abéliens sur $S$. Nous avons alors une suite exacte de groupes abéliens, déduite de la suite spectrale locale-globale pour les Ext (voir \cite{gro4}, exposé V, proposition 6.3, 3))
\begin{align*}
H^1(S,\Hom_S(F,F')) \rightarrow \ext(F,F') \rightarrow \Ext_S(F,F')(S) \rightarrow H^2(S,\Hom_S(F,F'))\,.
\end{align*}
En particulier, si $\Hom_S(F,F')$ est nul, alors $\Ext_S(F,F')(S)\simeq\ext(F,F')$. Par suite, le raisonnement étant encore valable après changement de base plat $T\rightarrow S$, on trouve que $\Ext_S(F,F')(T)\simeq \ext(F|_T,F'|_T)$ fonctoriellement en $T$. Enfin le lemme $\ref{null}$ nous permet d'appliquer ce raisonnement à la situation présente. 
\end{proof}

L'application $\delta$ peut être explicitée : soit $\Gamma\in \ext(A,\gm)$. Alors on peut associer à $\Gamma$, via l'isomorphisme $\ext(A,\gm)\simeq \Ext_S(A,\gm)(S)$, un élément de $H^1(S,G^D)$, que nous noterons $\delta(\Gamma)$ par abus de langage. Le lemme \ref{ext} permet de décrire $\delta(\Gamma)$ comme étant le faisceau
$$[\phi^*]^{-1}(\Gamma):T\longmapsto \{\Theta\in \ext(B_T,\gmt)\;|\;\phi^*\Theta=\Gamma_T\}\,.$$
En résumé, on peut caractériser $\delta(\Gamma)$ comme étant le faisceau des extensions $\Theta$ de $B$ par $\gm$ telles que $\phi^*\Theta=\Gamma$, ce qui justifie l'écriture $\delta(\Gamma)=[\phi^*]^{-1}(\Gamma)$.

\subsection{L'homomorphisme de Waterhouse}
Nous rappelons ici une construction, due à Waterhouse (voir \cite{w}, section 2), qui permet d'associer à toute extension $\Omega$ de $G$ par $\gm$ un $G^D$-torseur.
Considérons la suite exacte
\[
\begin{CD}
0 @>>> \gm @>>> \Omega @>>> G @>>> 0\,.\\
\end{CD}
\]
Elle donne lieu, par application du foncteur $\Hom_S(G,-)$, à une suite exacte
\[
\begin{CD}
0 @>>> G^D @>>> \Hom_S(G,\Omega) @>>> \Hom_S(G,G) @>>> 0\\
\end{CD}
\]
(rappelons que $\Ext_S(G,\gm)=0$ d'après \cite{gro7}, exposé VIII, 3.3.1). Par application du foncteur des sections, on obtient un morphisme $\underline\delta:\Homs(G,G)\rightarrow H^1(S,G^D)$. On note alors $\rho(\Omega)$ le $G^D$-torseur $\underline\delta({\rm id})$. Autrement dit, $\rho(\Omega)$ est le faisceau des sections $s:G\rightarrow \Omega$, au sens de la théorie des extensions. On définit ainsi un morphisme $\rho:\ext(G,\gm)\rightarrow H^1(S,G^D)$, et Waterhouse a montré que $\rho$ est un isomorphisme (voir \cite{w}, Theorem $2'$). Ainsi, en composant les morphismes suivants (où $l$ est le morphisme naturel)
\[
\begin{CD}
H^1(S,G^D) @>\rho^{-1}>> \ext(G,\gm) @>l>> \pic(G) \\
\end{CD}
\]
on obtient l'homomorphisme $\pi$ de Waterhouse, qui admet une interprétation galoisienne, comme nous l'avons précisé dans l'introduction.

\begin{lem}
\label{diagref}
Soit $i:G\rightarrow A$ l'inclusion. Alors le diagramme suivant :
\[
\begin{CD}
\pic(A) @<<< \ext(A,\gm) @>\delta>> H^1(S,G^D)\\
@VVV @Vi^*VV @| \\
\pic(G) @<l<< \ext(G,\gm) @>\rho>> H^1(S,G^D) \\
\end{CD}
\]
est commutatif.
\end{lem}

\begin{proof}
La commutativité du carré de gauche est claire par fonctorialité. Soit $\Gamma$ dans $\ext(A,\gm)$, on veut montrer que les torseurs $\delta(\Gamma)$ et $\rho(i^*\Gamma)$ sont isomorphes. D'après ce qui précède (lemme \ref{ext}), on peut transposer le problème dans la catégorie des groupes abéliens. Dans ce cadre, il est bien connu qu'étant donnés un groupe $\breve{\rm Z}$, une suite exacte (de groupes abéliens)
\[
\begin{CD}
0 @>>> \breve{\rm G} @>i>> \breve{\rm A} @>\phi>> \breve{\rm B} @>>> 0\,,\\
\end{CD}
\]
et une extension $\Gamma$ de $\breve{\rm A}$ par $\breve{\rm Z}$, alors il y a bijection entre les sections de $i^*\Gamma$ et les extensions $\Theta$ de $\breve{\rm B}$ par $\breve{\rm Z}$ telles que $\phi^*\Theta=\Gamma$. Le résultat découle alors des descriptions de $\delta$ et de $\rho$ que nous avons données précédemment.
\end{proof}

\subsection{Définition de l'homomorphisme}
\label{defpsi}
Nous sommes maintenant en mesure de généraliser la construction de Taylor. Reprenons les notations et hypothèses introduites dans le paragraphe \ref{bibi}.

On déduit du lemme \ref{diagref} le diagramme commutatif suivant :
\[
\begin{CD}
A'(S) @>\gamma>> \ext(A,\gm) @>\delta>> H^1(S,G^D) \\
@. @Vl^1 VV @VV\pi V \\
@. \pic(A) @>>> \pic(G) \\
\end{CD}
\]
et on définit $\psi:A'(S)\rightarrow\pic(G)$ comme étant le composé de ces morphismes. Dans le cas où $A$ est un $S$-schéma abélien, alors $A'$ est le schéma abélien dual de $A$, la flèche $\gamma$ est un isomorphisme, et on retrouve le {\em class-invariant homomorphism}, dont la première construction est due à M. J. Taylor (voir \cite{t2}).

\begin{rmq}
\label{geom}
Notons $\mathcal{D}:A'(S)\rightarrow\pic(A)$ la composée $l^1\circ\gamma$. Alors le diagramme précédent montre que, pour tout $p\in A'(S)$, $\psi(p)$ est la restriction de $\mathcal{D}(p)$ à $G$. Ceci généralise la << description géométrique >> de $\psi$ (obtenue par Agboola \cite{a1} dans le cas où $A$ est un schéma abélien).

Nous pouvons donner une autre écriture de $\psi$ : soit $p\in A'(S)$, alors on a l'égalité
$$\mathcal{D}(p)=l^1(({\rm id}_A\times p)^*(W))=({\rm id}_A\times p)^*(t(W))\,.$$
On en déduit que $\psi(p)=(i\times p)^*(t(W))$, où $i:G\rightarrow A$ est l'inclusion.
\end{rmq}

\begin{prop}
\label{oo}
Soit $\ord(G)$ l'ordre de $G$. Alors $\ord(G).\psi=0$. En particulier, $\psi$ s'annule sur les points de torsion d'ordre premier à $\ord(G)$.
\end{prop}

\begin{proof}
Le groupe $G$ étant fini et plat, la multiplication par l'entier $\ord(G)$ est l'application nulle dans $G$, donc est également l'application nulle dans le groupe $\ext(G,\gm)$. Or $\psi$ se factorise à travers ce dernier, d'où le résultat.
\end{proof}

\begin{rmq}
\label{imp}
Soit $p\in A'(S)$ un point de $m$-torsion. On peut écrire $m=nN$ où $N$ est premier à $\ord(G)$, et l'ensemble des facteurs premiers de $n$ est un sous-ensemble de l'ensemble des facteurs premiers de $\ord(G)$. Alors $Np$ est un point de $n$-torsion, et il découle de la proposition \ref{oo} que la nullité de $\psi(Np)$ implique celle de $\psi(p)$.

Quitte à changer $n$ en un multiple, dont les facteurs premiers sont ceux de $G$, on peut supposer que $G$ est un sous-groupe de $A[n]$. Sachant que $Np$ se factorise à travers $A'[n]$, on peut alors écrire (cf. la remarque \ref{geom})
$$\psi(Np)=(i\times Np)^*(t(W)|_{A[n]\times_S A'[n]})\,.$$
Ainsi, pour montrer que $\psi(p)$ est nul, il suffit de montrer que la restriction de $t(W)$ à $A[n]\times_S A'[n]$ est triviale.
\end{rmq}

\subsection{Constructions d'exemples}
\label{explicit}
Pour engendrer un exemple, on doit trouver un sous-$S$-schéma en groupes fini et plat $G$ de $A$. Une idée naturelle est de regarder ce qui arrive si les points de $N$-torsion de $A_{\eta}$ sont rationnels sur $K$. Fixons une clôture algébrique $\overline{K}$ de $K$. Nous avons le résultat suivant (voir \cite{mb}, chap. IV, cor. 8.2) :

\begin{prop}
\label{moret}
Soient $N$ un entier naturel, $F_1$ le corps $K(A_{\eta}(\overline{K})[N])$, et $v:S_1\rightarrow S$ le normalisé de $S$ dans $F_1$. Alors $S_1$ est intègre normal, et $v$ est fini surjectif de rang divisant le cardinal de $\mathbf{GL}_{2g}(\mathbb{Z}/N\mathbb{Z})$, où $g$ est la dimension de $A_{\eta}$. En outre, il existe un $S_1$-schéma en groupes semi-stable $A_1^{\sharp}$, contenant $A_1:=A\times_S S_1$ comme sous-groupe ouvert et ayant même fibre générique que lui, tel que le noyau de la multiplication par $N$ dans $A_1^{\sharp}$ soit fini et plat sur $S_1$. Le groupe $A_1^{\sharp}(S_1)[N]$ est isomorphe à $A_{\eta}(\overline{K})[N]$.
\end{prop}

Notre énoncé diffère légèrement de celui de \cite{mb} ; les détails que nous avons rajoutés découlent naturellement de la démonstration de ce dernier. Remarquons également que l'hypothèse d'excellence de $S$ est utilisée dans ladite démonstration. Pour notre part, c'est le seul endroit où nous en ferons usage.

De façon analogue, tout point de torsion pris dans $A(S)$ engendre un sous-groupe fini et plat de $A$, comme l'affirme la proposition suivante :

\begin{prop}
\label{satge}
Soit $x\in A(S)$ un point d'ordre fini $m$. Alors $x$ définit un morphisme $(\mathbb{Z}/m\mathbb{Z})_S\rightarrow A$ dont l'image est un sous-schéma en groupes fini et plat de $A$.
\end{prop}

\begin{proof}
Soit $f:(\mathbb{Z}/m\mathbb{Z})_S\rightarrow A$ un morphisme. Alors l'image de $f$ est un sous-groupe quasi-fini et plat de $A$. De plus, d'après \cite{ega2}, corollaire 5.4.3 (ii), l'image de $f$ est propre (car $A$ est séparé sur $S$, et $(\mathbb{Z}/m\mathbb{Z})_S$ est propre sur $S$), donc l'image de $f$ est finie d'après le Main Theorem de Zariski (quasi-fini et propre impliquent fini).
\end{proof}

Pour les courbes elliptiques, nous avons le critère suivant (voir \cite{maz1}, prop. 9.1) :

\begin{prop}
\label{mazur}
Supposons que $K$ soit un corps de nombres, et que $S=Spec(\mathcal{O}_K)$ soit le spectre de l'anneau des entiers de $K$. Soit $E\rightarrow S$ le modèle de Néron d'une courbe elliptique à réduction semi-stable définie sur $K$, de discriminant minimal $(\Delta)$, et soit $N$ un nombre premier. Alors le groupe $E[N]$ est fini et plat sur $S$ si et seulement si, pour toute place $\mathfrak{p}$ de $K$, $v_{\mathfrak{p}}(\Delta)$ est un multiple de $N$.
\end{prop}


\section{Cas des courbes elliptiques}

Dans cette section, nous donnons plusieurs applications de notre construction, dans lesquelles les variétés abéliennes en jeu sont des courbes elliptiques semi-stables.

\subsection{Torseurs rigidifiés}
\label{rigidifications}
Soit $f:X\rightarrow S$ un $S$-schéma en groupes. Soit $0_X$ la section unité de $X$, et soit $\mathcal{L}$ un $\gm$-torseur sur $X$. Une {\em rigidification} de $\mathcal{L}$ est la donnée d'une section du $\gm$-torseur $0_X^*\mathcal{L}$ sur $S$. On note $\Torsrig(X,\gm)$ la catégorie des $\gm$-torseurs rigidifiés sur $X$. C'est une catégorie de Picard strictement commutative (au sens de \cite{gro4}, exposé XVIII, 1.4.2), la loi étant le produit tensoriel. On note $\pic_r(X)$ le groupe des classes d'isomorphie d'objets de $\Torsrig(X,\gm)$.

Si $\mathcal{N}$ est un $\gm$-torseur sur $X$, on note
$$\mathbf{r}(\mathcal{N}):=\mathcal{N}\otimes (f^*0^*\mathcal{N})^{-1}$$
le $\gm$-torseur rigidifié associé à $\mathcal{N}$ (muni de sa rigidification naturelle).

Soient $g:Y\rightarrow S$ un autre $S$-schéma en groupes, $0_Y$ la section unité de $Y$, et $\mathcal{L}$ un $\gm$-torseur sur $X\times_S Y$. Une {\em birigidification} de $\mathcal{L}$ est la donnée de deux sections : une pour le torseur $(0_X\times {\rm id}_Y)^*\mathcal{L}$ sur $Y$ et une autre pour le torseur $({\rm id}_X\times 0_Y)^*\mathcal{L}$ sur $X$, qui soient compatibles entre elles (c'est-à-dire qui induisent la même section pour le torseur $(0_X\times 0_Y)^*\mathcal{L}$ sur $S$). On note $\Torsb(X,Y;\gm)$ la catégorie des $\gm$-torseurs birigidifiés sur $X\times Y$.

Si $\mathcal{N}$ est un $\gm$-torseur sur $X\times_S Y$, on note
$$\mathbf{bir}(\mathcal{N}):=\mathcal{N}\otimes(((0_X\circ f)\times {\rm id}_Y)^*\mathcal{N})^{-1}
\otimes(({\rm id}_X\times (0_Y\circ g))^*\mathcal{N})^{-1}\otimes (f\times g)^*(0_X\times 0_Y)^*\mathcal{N}$$
le $\gm$-torseur birigidifié associé à $\mathcal{N}$ (muni de sa birigidification naturelle).

\subsection{Démonstration du théorème \ref{jg}}
D'après la remarque \ref{imp}, pour montrer le théorème \ref{jg}, il suffit de montrer le théorème suivant :

\begin{thm}
\label{ppal}
Supposons que $A_{\eta}$ soit une courbe elliptique, et soit $n$ un entier premier à $6$. Alors le $\gm$-torseur $t(W)|_{A[n]\times_S A'[n]}\in \pic(A[n]\times_S A'[n])$ est trivial.
\end{thm}

Pour démontrer ce résultat, nous allons utiliser une méthode analogue à celle de Pappas dans \cite{p1}.

\begin{lem}
\label{lemme1}
Le $\gm$-torseur $(t(W)|_{A[n]\times_S A'[n]})^{\otimes n}$ est trivial.
\end{lem}

\begin{proof}
Il suffit de montrer que le groupe $\biext(A[n],A'[n];\gm)$ est tué par $n$. D'après (\cite{gro7}, exposé VIII, 1.1.2), on a un isomorphisme canonique :
$$\biext(A[n],A'[n];\gm)\simeq\ext(A[n],\mathbf{R}\Hom(A'[n],\gm))\,.$$
Par suite, la multiplication par $n$ étant nulle dans $A[n]$, elle est également nulle dans le groupe $\biext(A[n],A'[n];\gm)$. Ce qu'on voulait.
\end{proof}

\begin{lem}
\label{lemme2}
Pour démontrer le théorème \ref{ppal}, on peut supposer que $A'$ est à fibres connexes, et que le groupe $A(S)[6]$ est isomorphe à $A_{\eta}(\overline{K})[6]$, où $\overline{K}$ désigne une clôture algébrique de $K$.
\end{lem}

\begin{proof}
Quitte à remplacer $W$ par sa biextension symétrique $W^s$, on peut supposer que $A'$ est à fibres connexes. Appliquons la proposition \ref{moret} au schéma $A$ en prenant $N=6$. Alors $v:S_1\rightarrow S$ est de rang $k$ divisant $2^53^2$ et il existe un $S_1$-schéma en groupes semi-stable $A_1^{\sharp}$, contenant $A_1:=A\times_S S_1$ comme sous-groupe ouvert et ayant même fibre générique que lui, tel que $A_1^{\sharp}(S_1)[6]$ soit isomorphe à $A_{\eta}(\overline{K})[6]$. Soit $A'_1=A'\times_S S_1$, alors $A'_1$ est à fibres connexes. On note $W_1^{\sharp}\in \biext(A_1^{\sharp},A'_1;\gm)$ l'unique biextension dont l'existence est énoncée dans la proposition \ref{wo}. Alors la restriction de $W_1^{\sharp}$ à $(A_1, A'_1)$ est égale à $W_1:=W\times_S S_1$, par unicité du prolongement.

Pour simplifier, on pose $X:=A[n]\times_S A'[n]$, donc $X\times_S S_1=A_1[n]\times_{S_1} A'_1[n]$.
D'après ce qui précède, il suffit de montrer que la trivialité du $\gm$-torseur
$$t(W_1^{\sharp})|_{X\times_S S_1}=t(W_1)|_{X\times_S S_1}= (t(W)|_X)_{S_1}$$
implique la trivialité du $\gm$-torseur $t(W)|_{X}$.

Le schéma $S$ étant normal, il existe un gros ouvert $V$ de $S$ tel que $S_1\rightarrow S$ soit plat au-dessus de $V$. Nous adoptons ici la terminologie de (\cite{mb}, chap. II, déf. 3.1) : un gros ouvert $V$ de $S$ est un ouvert tel que, pour tout $x\in S\backslash V$, l'anneau local $\mathcal{O}_{S,x}$ soit de profondeur $\geq 2$. Soit $V_1:=S_1\times_S V$, alors $V_1\rightarrow V$ est fini localement libre de rang $k$.

Comme $X$ est un $S$-schéma plat, $X_V$ est un gros ouvert de $X$. Par suite, la flèche $\pic(X)\rightarrow \pic(X_V)$ est injective, d'après \cite{ega4}, 21.13.3 et 21.13.4. Nous avons alors un diagramme commutatif
$$
\begin{CD}
\pic(X)[n] @>>> \pic(X\times_S S_1) \\
@VVV @VVV \\
\pic(X_V)[n] @>>> \pic(X_V\times_V V_1) \\
\end{CD}
$$
où la flèche de gauche est injective. En se servant de la norme associée au morphisme $X_V\times_V V_1\rightarrow X_V$ qui est fini localement libre de rang $k$, ainsi que du fait que $k$ divise $2^53^2$, donc est premier à $n$, on montre que le flèche du bas est injective. On en déduit que la flèche du haut est injective. Or, d'après le lemme \ref{lemme1}, $t(W)|_X$ est annulé par $n$. On en déduit le résultat voulu.
\end{proof}

\begin{ntt}
Soit $X$ un $S$-schéma. Pour tout sous-schéma fermé $z:Z\subseteq X$ on note $I(Z)$ le faisceau d'idéaux définissant $Z$. Lorsque $I(Z)$ est inversible on note $I^{-1}(Z)$ son inverse. De plus, si $p:S\rightarrow X$ est un $S$-point de $X$, on note $\{p\}$ l'image de $p$ dans $X$, qui est sous-schéma fermé de $X$. Enfin, si $X$ est un $S$-schéma en groupes commutatif, on note $Z+p$ l'image du morphisme
$$
\begin{CD}
Z\times_S S @>z\times p>> X\times_S X @>+>> X
\end{CD}
$$
où $+$ désigne la loi de $X$. On vérifie que $Z+p$ est un sous-schéma fermé de $X$.
\end{ntt}

\begin{lem}
\label{lemme3}
Soit $E$ un $S$-schéma en groupes semi-stable dont la fibre générique $E_{\eta}$ est une courbe elliptique, et soit $\ez$ la composante neutre de $E$. On pose $A=E$ et $A'=\ez$, alors les hypothèses du paragraphe \ref{bibi} sont satisfaites, par auto-dualité de $E_{\eta}$.

D'autre part, soit $\Delta$ la diagonale de $E\times_S E$, et soit $(q_1,q_2)\in E(S)^2$. On note $\Pc(q_1,q_2)$ le $\gm$-torseur birigidifié sur $E\times_S E$ défini par
$$
\Pc(q_1,q_2)=\mathbf{bir}\left(I^{-1}(\Delta+(q_1,q_2))\otimes I(E\times\{q_2\})\otimes I(\{q_1\}\times E)\right)\,.
$$
Alors la biextension $W$ sur $(E,\ez)$ satisfait : $t(W)=\Pc(q_1,q_2)|_{E\times_S \ez}$.
\end{lem}

\begin{proof}
D'après la théorie classique des courbes elliptiques, le fibré de Poincaré $\mathcal{P}_{\eta}=t(W_{\eta})$ sur $E_{\eta}\times_K E_{\eta}$ s'écrit
$$t(W_{\eta})=I^{-1}(\Delta_{\eta})\otimes I(E_{\eta}\times\{0\})\otimes I(\{0\}\times E_{\eta})$$
où $\Delta_{\eta}$ est la diagonale de $E_{\eta}\times_K E_{\eta}$. De plus, nous avons l'égalité
$$
\Pc(q_1,q_2)_{\eta}=\mathbf{bir}\left(I^{-1}(\Delta_{\eta}+(q_1,q_2))\otimes I(E_{\eta}\times\{q_2\})\otimes I(\{q_1\}\times E_{\eta})\right)=\mathcal{P}_{\eta}\,.
$$
Il en résulte que $\Pc(q_1,q_2)$ est un $\gm$-torseur birigidifié sur $E\times_S E$ qui prolonge $t(W_{\eta})$. D'autre part le foncteur
$$\Biext(E,\ez;\gm)\longrightarrow \Torsb(E,\ez;\gm)$$
de la catégorie des biextensions de $(E,\ez)$ par $\gm$ dans la catégorie des $\gm$-torseurs birigidifiés sur $E\times_S \ez$, est pleinement fidèle (voir \cite{gro7}, exposé VIII, prop. 7.4, b)). De plus, $\ez$ est à fibres connexes, et $S$ est intègre normal, donc ({\em loc. cit.}) un objet $Z$ du second membre appartient à l'image (essentielle) du foncteur si et seulement si $Z_{\eta}$ provient d'une biextension de $(E_{\eta},E_{\eta})$ par $\gm$.

Soit $\Pc(q_1,q_2)|_{E\times_S\ez}\in\Torsb(E,\ez;\gm)$ la restriction de $\Pc(q_1,q_2)$ à $E\times_S\ez$. D'après ce qui précède, $\Pc(q_1,q_2)|_{E\times_S\ez}$ provient d'une biextension de $(E,\ez)$ par $\gm$, laquelle prolonge $W_{\eta}$, donc est égale à $W$ en vertu de la proposition \ref{wo}. Au final, $t(W)$ est égal à $\Pc(q_1,q_2)|_{E\times_S\ez}$ en tant que $\gm$-torseur birigidifié sur $E\times_S\ez$.
\end{proof}

\begin{proof}[Démonstration du théorème \ref{ppal}]
On suppose ici que $A_{\eta}$ est une courbe elliptique, et l'on note $A=E$ pour le confort du lecteur. D'après le lemme \ref{lemme2}, on peut supposer que $A'$ est à fibres connexes. Par auto-dualité des courbes elliptiques, la fibre générique de $A'$ est isomorphe à $E_{\eta}$. Or il existe (au plus) un unique prolongement de $E_{\eta}$ en un $S$-schéma en groupes semi-stable à fibres connexes. Par suite, $A'$ est isomorphe à la composante neutre $\ez$ de $E$. On retrouve la situation du lemme \ref{lemme3}.

D'après le lemme \ref{lemme3}, pour montrer le théorème \ref{ppal}, il suffit de trouver $(q_1,q_2)\in E(S)^2$ tel que la restriction de $\Pc(q_1,q_2)$ à $E[n]\times_S E[n]$ soit triviale. D'après le lemme \ref{lemme2}, on peut supposer que $E(S)[6]$ est isomorphe à $E_{\eta}(\overline{K})[6]$.

Supposons que la caractéristique de $K$ soit première à $6$. Alors $E_{\eta}(\overline{K})[6]$ est isomorphe à $(\mathbb{Z}/6\mathbb{Z})^2$. On peut donc choisir deux points $q_1$ et $q_2$ d'ordre $6$ dans $E(S)$ tels que $q_1-q_2$ soit d'ordre $6$.

Montrons, en raisonnant fibre par fibre, que $\{q_1\}\times E$ est disjoint de $E[n]\times_S E[n]$. Soit $x$ un point de $S$, de corps résiduel $k(x)$ ; alors l'application de réduction $E(S)\rightarrow E(k(x))$ est injective sur les points d'ordre premier à la caractéristique de $k(x)$. Ainsi, $q_1$ peut se réduire selon les cas en un point d'ordre $2$, $3$ ou $6$. L'entier $n$ étant premier à $6$, on en déduit que $\{q_1\}$ ne rencontre pas $E[n]$, d'où le résultat.

Un raisonnement analogue montre que $\Delta+(q_1,q_2)$ et $E\times\{q_2\}$ sont égalements disjoints du sous-schéma $E[n]\times_S E[n]$. On en déduit aussiôt que la restriction de $\Pc(q_1,q_2)$ à $E[n]\times_S E[n]$ est triviale.

Supposons que la caractéristique de $K$ soit égale à $2$. Alors $E_{\eta}(\overline{K})[3]$ est isomorphe à $(\mathbb{Z}/3\mathbb{Z})^2$. On peut donc choisir deux points $q_1$ et $q_2$ d'ordre $3$ dans $E(S)$ tels que $q_1-q_2$ soit d'ordre $3$. En outre, $S$ étant intègre, les caractéristiques résiduelles des points de $S$ sont toutes égales à $2$, donc l'argument utilisé précédemment permet de montrer que la restriction de $\Pc(q_1,q_2)$ à $E[n]\times_S E[n]$ est triviale.

Enfin, si la caractéristique de $K$ est égale à $3$, on effectue un raisonnement analogue en considérant les points d'ordre $2$.
\end{proof}

\begin{rmq}
Dans les hypothèses du paragraphe \ref{bibi}, $A$ ou $A'$ est à fibres connexes. Si l'on se contente de supposer que les groupes de composantes de $A$ et $A'$ sont orthogonaux sous l'accouplement de monodromie, alors $W_{\eta}$ admet à nouveau un (unique) prolongement en une biextension $W$ de $(A,A')$ par $\gm$. Nous ignorons si le théorème \ref{ppal} est encore vrai dans ce cadre plus général. Signalons simplement que le lemme \ref{lemme3} ne se généralise pas si l'on supprime l'hypothèse de connexité.
\end{rmq}

\subsection{Un exemple elliptique}
Donnons un exemple d'application du théorème \ref{jg}. Supposons que $K$ soit un corps de nombres, que $S$ soit le spectre de l'anneau des entiers de $K$, et que $A=E$ soit le modèle de Néron sur $S$ de la courbe elliptique $E_{\eta}:=X_0(11)$ (notée A1(B) par Cremona \cite{cr}) définie sur $K$ par l'équation
$$y^2+y=x^3-x^2-10x-20\,.$$
Le discriminant de cette courbe est $-11^5$, et $E$ est un $S$-schéma en groupes semi-stable. Dans le cas présent, on constate (voir \cite{maz1}, p. 258) que $\ez$ contient un sous-groupe isomorphe à $\boldsymbol{\mu}_{5/S}$.
Soit à présent
$$\psi:E(S)\rightarrow \pic(\boldsymbol{\mu}_{5/S})$$
l'homomorphisme de classes correspondant à $G=\boldsymbol{\mu}_{5/S}\subseteq \ez$ et $A'=E$. Le théorème \ref{jg} affirme que $\psi$ s'annule sur les points de torsion.

Du point de vue de la structure galoisienne, ce résultat peut s'interpréter de la façon suivante : à tout point $p\in E(S)_{\rm Tors}$ nous avons associé un $(\mathbb{Z}/5\mathbb{Z})_S$-torseur qui a une structure galoisienne triviale. En fait, si $p$ est d'ordre premier à $5$, le torseur associé est lui-même trivial. Les exemples intéressants sont donc fournis par des points $p$ d'ordre une puissance de $5$ : on pourra prendre $K=\mathbb{Q}(E_{\eta}[5^n])$ avec $n\geq 1$ dans cet exemple.

D'autre part, la proposition \ref{satge} permet de construire d'autres exemples de sous-groupes de $\ez$, lesquels engendrent à leur tour d'autres torseurs.

\subsection{Dualité}
Nous allons donner une interprétation de notre résultat en termes de fibrés en droites. Tout d'abord, on note $\pic_r^0(A_{\eta})$ le sous-groupe de $\pic_r(A_{\eta})$ constitué par les classes des torseurs qui sont algébriquement équivalents à $0$.

D'autre part, on note $\pic_r^0(A)$ l'image réciproque de $\pic_r^0(A_{\eta})$ par le morphisme de restriction à la fibre générique $\mathcal{R}:\pic_r(A)\rightarrow\pic_r(A_{\eta})$.

La théorie habituelle de la dualité pour les variétés abéliennes affirme que l'application
\begin{equation}
\label{fld}
\begin{split}
A_{\eta}^t(K)\;\longrightarrow & \;\pic_r^0(A_{\eta})\\
(p:Spec(K)\rightarrow A_{\eta}^t)\;\longmapsto & \;({\rm id}_{A_{\eta}}\times p)^*(\mathcal{P}_{\eta})\\
\end{split}
\end{equation}
est un isomorphisme de groupes.

Nous allons voir que, sous certaines hypothèses, la biextension $W$ permet d'établir une dualité du même type entre $A$ et $A'$.

\begin{prop}
\label{linebundle}
Supposons que $A$ soit à fibres connexes, et que $A'$ soit néronien (déf. \ref{ssneron}). Alors l'application
\begin{equation}
\label{fldx}
\begin{split}
A'(S)\;\longrightarrow & \;\pic_r^0(A)\\
(p:S\rightarrow A')\;\longmapsto & \;\mathcal{D}(p)=l^1(({\rm id}_A \times p)^*(W))\\
\end{split}
\end{equation}
est un isomorphisme de groupes.
\end{prop}

\begin{proof}
Rappelons que $\mathcal{D}=l^1\circ\gamma$ (voir le paragraphe \ref{bibi}). D'après la remarque \ref{gamma}, les hypothèses que nous avons faites impliquent que l'application $\gamma$ est un isomorphisme de groupes. D'autre part, le foncteur
$$L^1:{\rm EXT}(A,\gm)\longrightarrow \Torsrig(A,\gm)$$
(où ${\rm EXT}(A,\gm)$ désigne la catégorie des extensions de $A$ par $\gm$) est pleinement fidèle (voir \cite{gro7}, exposé VIII, prop. 7.4, a)). De plus, $A$ est à fibres connexes, et $S$ est intègre normal, donc ({\em loc. cit.}) un objet $Z$ du second membre appartient à l'image (essentielle) du foncteur $L^1$ si et seulement si $Z_{\eta}$ provient d'une extension de $A_{\eta}$ par $\gmk$.

Considérons alors le diagramme commutatif suivant
\[
\begin{CD}
A'(S) @>\gamma>> \ext(A,\gm) @>l^1>> \pic_r(A) \\
@| @VVV @VV\mathcal{R}V \\
A_{\eta}^t(K) @>\sim>> \ext(A_{\eta},\gmk) @>>> \pic_r(A_{\eta}) \\
\end{CD}
\]
dans lequel $l^1$ est l'application induite par le foncteur $L^1$, et $\mathcal{R}$ est l'application de restriction à la fibre générique. Alors, d'après ce qui précède, $l^1$ est injective, et son image est égale à l'image réciproque par $\mathcal{R}$ de l'image de l'application (\ref{fld}). Or cette image est égale à $\pic_r^0(A_{\eta})$, d'où le résultat.
\end{proof}

En combinant la proposition \ref{linebundle} et le théorème \ref{ppal}, nous obtenons le résultat suivant, qui généralise le Theorem A de Pappas \cite{p1} :

\begin{cor}
\label{quasifini}
Supposons que $A$ soit à fibres connexes, et que $A_{\eta}$ soit une courbe elliptique. Soit $m$ un entier premier à $6$. Alors, pour tout $\mathcal{L}$ dans $\pic_r^0(A)_{\rm Tors}$, la restriction de $\mathcal{L}$ à $A[m]$ est triviale.
\end{cor}

\begin{rmq}
En général, le schéma en groupes $A[m]$ est quasi-fini et plat sur $S$, mais n'est pas nécessairement fini. Cependant, dans le cas particulier où $S$ est le spectre d'un anneau de Dedekind, le corollaire \ref{quasifini} admet encore une interprétation en termes de structure galoisienne de torseurs. Ceci fait l'objet d'un travail en cours de préparation.
\end{rmq}

\subsection{Application aux caractéristiques d'Euler}
Comme l'a montré Pappas dans \cite{p2}, l'homomorphisme $\psi$ admet une interprétation en termes de caractéristiques d'Euler. Ainsi deux problèmes de structure galoisienne, {\em a priori} distincts, se retrouvent liés --- la connexion étant établie par le biais de la théorie des biextensions.

Plus précisément, dans les paragraphes $3.b$ et $3.c$ de \cite{p2}, l'auteur étudie les caractéristiques d'Euler dans un cadre général. Puis, sous l'hypothèse de bonne réduction, il établit un lien entre caractéristiques d'Euler et valeurs de l'homomorphisme $\psi$. Dans ce cadre, son Theorem 6.4 (voir la section 6 de \cite{p2}) découle de l'annulation de l'homomorphisme $\psi$.
On vérifie que les arguments de Pappas s'appliquent {\em mutatis mutandis} à notre situation. On peut alors déduire du théorème \ref{jg} le résultat suivant, généralisant le Theorem 6.4 de Pappas :

\begin{cor}
Soit $S$ le spectre d'un anneau de Dedekind, et soit $G$ un $S$-schéma en groupes commutatif, fini et plat.

Soit $Y\rightarrow S$ un modèle minimal, régulier et projectif d'une courbe elliptique à réduction semi-stable sur $K$, et soit $X\rightarrow Y$ un $G$-torseur. Alors
$$pgcd(12m,2^7\,3^2)\cdot\chi_R^P(\mathcal{O}_X)=0$$
où $m$ est l'ordre de $G$, et où $\chi_R^P$ désigne la caractéristique d'Euler projective équivariante. De plus, si $m$ est premier à $6$ et si $R$ est principal, alors $\chi_R^P(\mathcal{O}_X)=0$.
\end{cor}

Pour une définition de la caractéristique d'Euler projective $\chi_R^P$, nous renvoyons le lecteur à l'article de Pappas (voir \cite{p2}, paragraphes $2.b$ et $2.c$).


\end{document}